\documentclass[12pt,titlepage]{article}

\usepackage{amssymb}
\usepackage{natbib}
\bibpunct{(}{)}{;}{a}{,}{,}
\usepackage{graphicx}
\usepackage{url}
\usepackage{verbatim} 

\usepackage[tbtags]{amsmath}
\usepackage[onehalfspacing]{setspace}
\usepackage{lineno} 
\usepackage[all]{xy} 
\usepackage{xspace}

\newcommand{\yobs}{y^\mathrm{obs}}
\newcommand{\curlyM}{\mathcal{M}}
\newcommand{\curlyY}{\mathcal{Y}}
\newcommand{\R}{\mathbb{R}}
\newcommand{\rd}{\mathrm{d}}
\newcommand{\ldef}{\mathrel{:=}\nolinebreak}
\newcommand{\that}{\hat{\theta}}

\newcommand{\thalf}{\tfrac{1}{2}}
\newcommand{\trans}{^{\scriptscriptstyle T}}

\newcommand{\Eqref}[1]{Eq.~\eqref{#1}}

\newcommand{\conv}[1]{\xrightarrow[\hphantom{1.25cm}]{\textrm{#1}}\nolinebreak
}
\newcommand{\Pconv}{\conv{P}} 
\newcommand{\pen}{\textsf{penalty}\xspace}
\newcommand{\pprime}{\text{$\pen^{\prime}$}}
\newcommand{\flex}{\textsf{flexibility}\xspace}

\DeclareMathOperator*{\argmax}{argmax}

\DeclareMathOperator{\Normal}{\mathcal{N}\!}

\DeclareMathOperator{\bigO}{O}

\begin{document}

\title{The exact form of the `Ockham factor'\\in model selection}

\author{Jonathan Rougier\thanks{%
Jonathan Rougier (corresponding author) is Honorary Professor, School
of Mathematics, University of Bristol (email:
\texttt{j.c.rougier@bristol.ac.uk}); Carey Priebe is
Professor, Department of Applied Mathematics and Statistics, Johns
Hopkins University (email: \texttt{cep@jhu.edu}).}\\
       School of Mathematics \\
       University of Bristol, UK
       \and
       Carey E. Priebe \\
       Department of Applied Mathematics and Statistics \\
       Johns Hopkins University, USA}

\date{ }


\maketitle


\begin{abstract}%

We explore the arguments for maximizing the `evidence' as an
algorithm for model selection.  We show, using a new definition of
model complexity which we term `flexibility', that maximizing the
evidence should appeal to both Bayesian and Frequentist statisticians.
This is due to flexibility's unique position in the exact decomposition
of log-evidence into log-fit minus flexibility.  In the Gaussian linear
model, flexibility is asymptotically equal to the Bayesian Information
Criterion (BIC) penalty, but we caution against using BIC in place
of flexibility for model selection.

\bigskip\noindent%
Keywords:
  Complexity, evidence, Ockham factor, flexibility, BIC penalty

\end{abstract}

\section{Introduction}
\label{sec:intro}

As Albert~Einstein almost said, ``Scientific models should be as
simple as possible, but no simpler''. \emph{Ockham's razor} is a
common term for this principle of simplicity. \citet{sober15} contains
a wide-ranging review, including what Einstein actually said, while
\citet{jefferys92} give an engaging and non-technical account of how
Ockham's razor emerges from a Bayesian treatment of model selection.

In philosophy, statistics, and machine learning, Ockham's razor
recommends that we favour less complex models where we can.  But why
this should be, and what we mean by `complex', are subtle issues.
In this paper, we describe one approach to implementing Ockham's
razor for model selection, based on a decomposition of the `evidence'
into a `fit' term and a complexity term which we call `flexibility'.
Our decomposition is exact for all models and all regularizers.
We do not advocate one particular approach to model selection,
but clarify the implications of and relations between various approaches,
describing how one particular approach---maximizing the
evidence---might be justified by both Bayesian and Frequentist
statisticians.

Section~\ref{sec:back} provides background, and
section~\ref{sec:arg} provides and discusses various arguments for the
use of the evidence in model selection.  Section~\ref{sec:candidate}
presents our exact decomposition of evidence into `fit plus
flexibility', and we justify flexibility as a measure of model
complexity.  Section~\ref{sec:GauLM} illustrates using the Gaussian
linear model, for which simple closed-form expressions are available.
The Bayesian Information Criterion (BIC) penalty is shown to be an
approximation to flexibility in this case, but we caution against
its use in model selection in general, and even specifically in those
cases where the BIC penalty and flexibility are asymptotically
equivalent.
Section~\ref{sec:summary} concludes with a summary.

\section{Background}
\label{sec:back}

Our starting-point is a set of observations $\yobs$.  A statistical model is proposed,
\begin{equation}
  f(y; \theta) , \qquad
  \begin{cases}
    y \in \curlyY \subset \R^n \\
    \theta \in \Theta \subset \R^d .
  \end{cases}
\end{equation}
The defining features of such a model are
\begin{equation}
  \forall \theta \in \Theta: \quad f(\cdot; \theta) \geq 0, \, \int_{\R^n} f(y; \theta) \, \rd y = 1 ,
\end{equation}
where the integral may be replaced by a sum if $\curlyY$ is countable.
In addition, this model is augmented with a regularizer $R : \Theta \to \R$.
The fitted value of the parameter is
\begin{equation}\label{eq:that}
  \that \ldef \argmax_{\theta \in \Theta} \, \big\{ \log f(\yobs; \theta) - R(\theta) \big\} .
\end{equation}
Thus the regularizer is functionally equivalent to a prior distribution
\begin{equation}\label{eq:prior}
  \pi(\theta) \propto e^{-R(\theta)} ,
\end{equation}
which we assume is proper, and in this case the fitted value $\that$
is identical to the Maximum A Posteriori (MAP) estimate of $\theta$.

As this outline makes clear, the estimated value $\that$ is a function
of both the statistical model and the regularizer.  To be precise, we would write `model and regularizer' everywhere below, but this would
be tedious; therefore we will write `model', but in every case where
we write `model' we mean `model and regularizer'.

Define the `evidence' of this model as
\begin{equation}
  E \ldef \int_{\Theta} f(\yobs; \theta) \, \pi(\theta)\, \rd \theta .
\end{equation}
This value is also referred to as the `integrated likelihood', or the
`marginal likelihood' \citep[sec.~5.3]{murphy12}.
In a Bayesian approach it has the additional interpretation of
the `predictive density' of the observations \citep[p.~385]{box80}.  
The normalizing constant in the prior distribution
\eqref{eq:prior} is not required to compute $\that$, but it \emph{is}
required to compute the evidence directly.  \citet{friel12}
provide a review of methods for estimating the evidence, covering a
literature which stretches back thirty years, while \citet{gutmann12}
present a recent and promising approach using logistic regression
\citep[see
also][sec.~14.2.4]{hastie09}.

Suppose that we would like to
select a single model from a set of models under consideration,
indexed by $i \in \curlyM$.  For example, $\curlyM$ might represent a
set of Gaussian linear models with different model matrices (see
Section~\ref{sec:GauLM}).  The claim we investigate in this paper is that maximizing the
evidence is a good way to select such a model.  That is, if $E_i$
is the evidence of model $i$, then 
\begin{equation}
  i^* \ldef \argmax_{i \in \curlyM} E_i
\end{equation}
is the best single model in $\curlyM$---or one of the best.

The notion that there is a best single model in $\curlyM$ is a bold
claim, because it does not distinguish between the possibly different
requirements of model-selection for inference, and model-selection
for prediction.  In inference, the intention would be to inspect the
fitted components of the model in order to make statements about
the underlying system; for example, to see whether a regression
coefficient is negative, roughly zero, or positive.  In this case,
it is helpful for interpretability to focus on a single model.
In prediction, the structure of the model is less important than
its predictive performance out-of-sample.  In this case, the use of
a single model would be a pragmatic simplification, especially for
operational systems.  So in both cases there is an apparent need for
a single model, but there is no \textit{a priori} reason to think
that one method will produce a best or nearly-best choice of model
in both cases.

\section{Arguments for the evidence}
\label{sec:arg}

For statisticians, the natural approach to choosing a single model
from a set $\curlyM$
is to use statistical decision theory, in which a loss function $L(i,
j)$ quantifies the negative consequences of choosing model $i$ were
model $j$ to be the true model.  While there will be 
applications where a particular loss function suggests itself, in
many cases the zero-one loss function is a pragmatic choice:
\begin{equation}
  L(i, j) = \begin{cases} 0 & i = j \\ 1 & \text{otherwise.} \end{cases}
\end{equation}
This loss function treats all wrong choices as equally-wrong.  There is
a useful generalization of this loss function, which includes
`undecided' among the actions, but we will not consider it further
here \citep[see][ch.~5]{murphy12}.

In a Bayesian approach,
the models in $\curlyM$ are given prior probabilities, $w_i > 0$.
The Bayes Rule for zero-one loss is to select the model with the
largest posterior probability, where the posterior probabilities are
proportional to $w_i \, E_i$.  Thus if $\curlyM$ is finite and
the prior
probabilities are uniform, then choosing $i^*$, the model with
the largest
evidence, is the Bayes Rule.  More generally, if the
evidence is concentrated relative to the prior probabilities,
then choosing $i^*$ is `hopefully' the Bayes Rule.  A Bayesian
statistician might reason: ``I don't want to spend too much effort
thinking about my prior probabilities, but I believe that they are
sufficiently
uniform over $\curlyM$ that that the model with the
highest evidence will be the Bayes Rule, or a good approximation.''
This reasoning goes back to the `Principle of Stable Estimation' of
L.J.~Savage \citep[see][]{edwards63}.

Clearly, this Bayesian argument is contingent on a willingness
to contemplate and possibly provide prior probabilities for the
models in $\curlyM$.  This may be anathema to some Frequentist
statisticians, who might otherwise accept the need for a loss function,
and tentatively accept the choice of the zero-one loss function.
But these statisticians will recognize the validity of the Complete
Class theorems, which in this context state that choosing the model
with the largest evidence cannot be inadmissible under zero-one loss,
because it is the Bayes Rule for a prior which assigns positive
probability to all models; see, e.g., \citet[ch.~8]{degroot70},
\citet[ch.~3]{schervish95}, or \citet[ch.~8]{berger85}.  Avoiding
inadmissible decision rules is a low bar but an important one.
So statisticians of all persuasions should at least contemplate
choosing the model with the highest evidence, on the basis of zero-one
loss, but this argument is suggestive, rather than compelling; and more
suggestive in the Bayesian approach than in the Frequentist approach.

The second argument is more Frequentist,
because it claims that $i^*$ is a model-selection algorithm with
attractive properties.  Making this case for
inference tends to be model-dependent and often requires asymptotic
arguments; we concur with \citet[p.~165]{lecam90}, whose Principle~7
states
``If you need to use asymptotic arguments, do not forget to let your
number of observations tend to infinity''; we touch on this principle
at the end of section~\ref{sec:GauLM}.  Instead, we
focus on the argument for choosing $i^*$ for prediction, using
a canonical example which also has much practical relevance:
polynomial regression.

In the regression
\[
  y_i = \beta_0 + \beta_1 \, x_i + \beta_2 \, x_i^2  + \cdots + \beta_p \, x_i^p + \text{error}_i
\]
there is a simple
index of models in terms of the degree $p$.  Under
standard fitting procedures such as Ordinary Least Squares, and its
Bayesian treatment with a vague prior on $\beta$, models with
higher degree will never fit worse, in terms of root mean squared
error, and will usually fit better.
Nevertheless, the wiggliness of high-degree polynomials
compromises their predictive accuracy, particularly
outside the convex hull of the observations.  This property is captured
by the following schematic:
\[
  \xymatrix{
    \txt{Low degree\\poor fit\\poor prediction} \ar[rrrr] &&&&
    \txt{High degree\\excellent fit\\poor prediction} \\
    & \txt{Medium degree\\good fit\\good prediction} \ar[u]
  }
\]
in which there is a `sweet spot' between low and high degree which
is optimal for prediction.  Most textbooks contain a detailed
discussion of this property; see, e.g., \citet[ch.~1]{bishop06}
or \citet[ch.~1]{murphy12}.  \citet{jefferys92} provides a wide-ranging
and non-technical discussion with historical examples.

The logic of polynomial regression would seem to apply much more
generally: models with the equivalent of high degree will fit better
but predict worse than models of medium degree.  This can be represented
as `overfitting': models of high degree fit into the 
noise, and thereby carry the noise into the prediction, which 
decreases accuracy.  Or, to put it another way, models need some
`stiffness' in order to resist noise.  The `degree' of polynomial
regression is replaced by the more ambiguous term `complexity',
and the schematic becomes 
\[
  \xymatrix{
    \txt{Simple model\\poor fit\\poor prediction} \ar[rrrr] &&&& \txt{Complex model\\excellent fit\\poor prediction} \\
    & \txt{Medium complexity\\good fit\\good prediction} \ar[u]
  }
\]
see, e.g., \citet[Fig.~2.11]{hastie09}.
It is only a short step from this schematic to an outline model-selection
criterion for prediction: to reward fit but to
penalize complexity, hoping to settle at or near the sweet-spot
between low and high complexity.

While there
are many forms that such a criterion might take, one form has emerged
as nearly ubiquitous, which is to maximize the scalar quantity
\begin{equation}\label{eq:criterion}
  \log f(\yobs; \that) - \pen
\end{equation}
where \pen is a measure of complexity that can depend on the model
and on the observations.  The first author to derive a criterion of
this form was \citet{akaike73}, whose penalty did not depend on the
data, and there have been many proposals since:
this is still an active field of research in statistics \citep[see,
e.g.,][ch.~7]{gelman14}.

\begin{figure}

\centering

\includegraphics[width=0.7\textwidth]{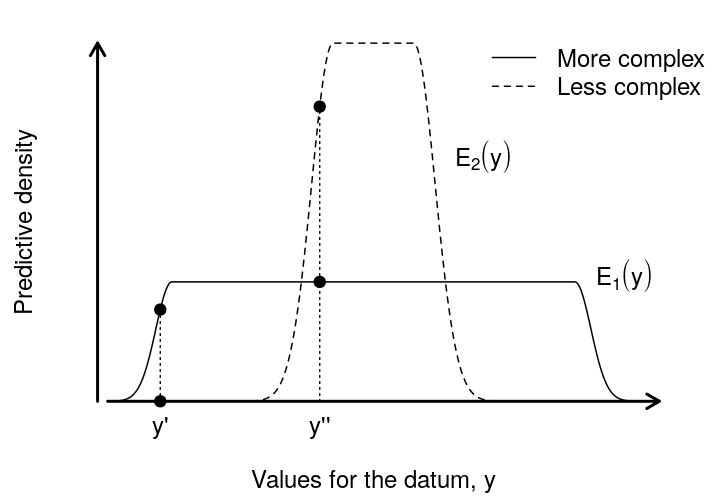}

\caption{Similar to \citet{mackay03}, Figure~28.3 (p.~344).
For visualization, the observation $y$ is assumed to be a scalar.  The more
complex Model~$1$ spreads its predictive density over a wider set of values
in the data-space $\curlyY$, and hence the evidence will favour the
less complex Model~$2$ in some regions of $\curlyY$.  In this case
the evidence
favours the more complex model at $\yobs = y'$, and the less complex
model at $\yobs = y''$.}

\label{fig:mackay}

\end{figure}

Now to return to our main topic, evidence.  David MacKay
\citep{mackay92,mackay03} argued that the log-evidence itself
has the form of \eqref{eq:criterion}, at least approximately.  
MacKay's argument had two strands.
First, there was `proof by picture', a compelling illustration that
the evidence will sometimes select less complex models over more
complex models, shown here as Figure~\ref{fig:mackay}.  This picture
embodies \citet{jefferys92}'s informal definition of complexity:
``the more complex hypothesis tends to spread the probability more
evenly among all the outcomes'' (p.~68).

But the second strand is less compelling.  MacKay applied a
first-order Laplace approximation to compute the evidence, and
showed that, under this approximation,
\begin{equation}\label{eq:laplace}
  \log E \approx \log f(\yobs; \that) - \pen
\end{equation}
where \pen has an explicit form in terms of the model and the
observations, but its form is not important to our argument.  We know,
as a matter of logic, that \pen in \eqref{eq:laplace} must behave something like a
model complexity penalty.  This follows because the first term on the
righthand side will tend to be larger for more complex models, and
hence if the evidence can sometimes be smaller for more complex
models, as Figure~\ref{fig:mackay} demonstrates, then \pen
must sometimes be larger for more complex models.

However, Laplace approximations, of first or higher order \citep[see,
e.g.,][sec.~3.5]{robert99} can perform poorly in practice, which is
why there has been such sustained interest in developing MCMC methods
for computing posterior expectations over the last three decades \citep[see, e.g.,][for
some history]{andrieu03}.
Their performance today is likely to be worse than ever, in our
current era of massively over-parameterized models and exotic
regularizers (e.g.\ to promote sparsity).  Therefore MacKay's
decomposition is suggestive, but does not provide a strong reason to
believe that the evidence on its own represents a fit term minus a
complexity term suitable for all applications.

This long discussion of the Frequentist argument for using
evidence to choose a single predictive model indicates once again that the
argument is suggestive but not compelling.  There have been a number
of conjectural steps: the generalization from polynomial degree to
model complexity, the idea that the `sweet spot' of optimal complexity
can be found by penalizing the maximized log-likelihood, and then
MacKay's construction which suggests that log-evidence has the
approximate form of penalized log-likelihood.  Therefore it is
gratifying that in the next section we can show that there is an exact
equality between log-evidence, fit, and a complexity penalty, which
holds in complete generality, and where the penalty has a simple and
intuitive form.

\section{`Flexibility'}
\label{sec:candidate}

Our approach uses a simple result that is a
reformulation of Bayes's Theorem:
\begin{equation}\label{eq:candidate}
  E = \frac{ f(\yobs; \theta_0) \, \pi(\theta_0) }{ \pi^*(\theta_0) }
\end{equation}
which holds for all $\theta_0 \in \Theta$
for which $\pi^*(\theta_0) > 0$, where $\pi^*$ is the posterior
distribution.  \citet[p.~1314]{chib95} refers to this equality as 
the \emph{basic marginal likelihood identity (BMI)}; it also goes by
the name \emph{Candidate's formula}, after \citet{besag89}.  

If we set $\theta_0 = \that$ from \eqref{eq:that}, then we immediately deduce
\begin{equation}
  \log E = \log f(\yobs; \that) - \log \frac{ \pi^*(\that) }{
  \pi(\that) } 
\end{equation}
from \eqref{eq:candidate}.
By our argument at the end of the previous section, the second term
on the righthand side must behave something like a model complexity
penalty.
To identify it explicitly, we give it the name
\begin{subequations}\label{eq:decomp}
\begin{equation}\label{eq:flex}
  \flex \ldef \log \frac{ \pi^*(\that) }{ \pi(\that) } ,
\end{equation}
which will be positive, except possibly in pathological situations where
there is a conflict between the prior distribution and the likelihood.
We have, under this definition,
\begin{equation}
  \log E = \log f(\yobs, \that) - \flex
\end{equation}
\end{subequations}
an exact decomposition of the evidence, which holds for all models.
This is the unique decomposition for which the `fit' term is $\log
f(\yobs; \that)$.  A different estimate for $\theta$, such as the
Generalized Method of Moments (GMM) estimate, would give a different
penalty term in the decomposition of the evidence.  However, given that
the regularizer is operationally equivalent to a prior distribution,
the penalized likelihood (or MAP) estimate $\that$ seems the most
natural value to use for $\theta_0$.

We contend that \flex is a reasonable way to measure model
complexity. A model will be `flexible' if it contains a large number
of degrees of freedom, in the engineering sense (e.g.\ represented
by $\dim \Theta$, which might be the number of basis functions),
and if its parameters are unconstrained in the prior distribution
(or regularizer).  A model will be `inflexible' (or `stiff')
either if it contains few degrees of freedom, or if its parameters
are constrained in the prior distribution, or both.  A flexible
model will often be able to
concentrate its posterior probability into a small data-determined
region of its parameter space, relative to its prior probability,
and hence its \flex will be high.  An inflexible model in the same
situation will be prevented from concentrating, lacking either the
capacity to shape itself to the data, or else having the capacity but
being prevented from doing so by the prior probability.  Thus its
\flex will be low.  Complex models will typically be flexible, in
this sense, and simple models will be inflexible.

Our exact decomposition of the log-evidence has changed neither
the Bayesian nor the Frequentist arguments for using the evidence for
model selection, although it has sharpened MacKay's
construction.  So does it have any practical value?  We think it does.  \Eqref{eq:decomp} shows that
\begin{equation}
  \log f(\yobs; \that) - \pen
  = \log E - \underbrace{(\pen - \flex)}_{\ldef \, \pprime} .
\end{equation}
Therefore any method which penalizes the maximum of the log-likelihood
with \pen is equivalent to a method which penalizes the log-evidence
with \pprime.  A decision \emph{not} to set $\pen = \flex$ is
equivalent to a decision to penalize the
log-evidence, with a penalty of an explicit but perhaps unanticipated
form.

Our decomposition presents this type of decision in a new light.
Suppose you choose to penalize log-likelihood with something other than
\flex.  Your client, or your client's auditor, might well ask
``What deficiency of the evidence as a criterion for model selection 
are you addressing with your choice of \pprime?''  You might answer,
``The evidence has no epistemic value to me: I regard \eqref{eq:decomp}
as a purely mathematical result, and can provide a rationale for
my \pen in terms of log-likelihood.''  The difficulty with this
position is that log-evidence clearly \emph{does} have epistemic value,
as demonstrated by the Bayesian and admissibility arguments.  Thus any
rationale for penalizing log-likelihood is more compelling if it can
be cross-referenced to log-evidence.

This viewpoint also operates in the other direction: a decision to
penalize log-evidence is equivalent to a decision to penalize
log-likelihood with something other than \flex.

There will always be applications in model selection where a strong case
can be made for a penalty which is not
\flex.  But we propose that \flex is
a reasonable default choice, which has some cross-party appeal.
Thus a Bayesian might say,
\begin{quote}
We selected the model which maximized the
evidence, which is similar and sometimes the same as selecting the
model which has the largest posterior probability, which itself is the Bayes Rule for zero-one loss.  But we note that
the same model would have been selected using a penalized log-likelihood
approach in which the model complexity penalty is `flexibility'.
\end{quote}
And a Frequentist might say
\begin{quote}
We selected the model which maximized
the log-likelihood penalized by the `flexibility' penalty for model
complexity; in other words, we maximized the evidence.  But we note
that the same model would have been selected 
under a Bayesian approach which maximized posterior probability with a uniform prior distribution, which is an admissible decision rule for zero-one loss.
\end{quote}

In situations where the choice of complexity penalty is not clear-cut,
choosing \flex seems to be a simple and defensible way forward.
The next section considers the case where the \flex has a
closed-form expression, and a simple asymptotic approximation.

\section{Illustration: Gaussian linear model}
\label{sec:GauLM}

Consider the Gaussian linear model
with model-matrix $G \in \R^{n \times d}$ and observation error
variance $\sigma^2$.  We use a quadratic regularizer parameterized by~$\lambda$:
\begin{subequations}
\begin{align}
  R(\theta) & = \thalf \, \lambda^2 \Vert \theta \Vert^2 \label{eq:reg} \\
\intertext{or, equivalently in terms of the prior distribution,}
  \pi(\theta) & = \Normal \big( 0_d, P^{-1})  \\
\intertext{where $P \ldef \lambda^2 I_d$ is the prior precision.  This implies that the posterior distribution has the form}
  \pi^*(\theta) & = \Normal \big( \that, (P^*)^{-1} \big) \\
\intertext{where $\that$ is the posterior expectation, as well as the MAP estimate, and $P^*$
is the posterior precision,}
  P^* & = \frac{1}{\sigma^2} G\trans G + P \\
  \that & = \frac{1}{\sigma^2} (P^*)^{-1} G\trans \yobs ;
\end{align}
\end{subequations}
these are both functions of $(\sigma, \lambda)$, which we treat as known.  Hence
\begin{equation}\label{eq:exact}
  \flex = \frac{1}{2} \log \left( \frac{ \det P^* }{ \det P } \right) + \frac{\lambda^2}{2} \Vert \that \Vert^2 .
\end{equation}
This is the exact result.  In the limit as $\lambda \to \infty$, $P^*
P^{-1} = I$ and $\that = \bigO(1/\lambda^2)$; thus $\flex
= \bigO(1 / \lambda^2)$, confirming that, asymptotically, \flex 
decreases to zero as the penalty on the quadratic regularizer
increases, although $\lambda$ might have to be huge to overwhelm the
$G\trans G$ term in $P^*$.

Now consider the effect of $n \to \infty$ when $d$ is fixed.
Under IID sampling,
\begin{subequations}
\begin{align}
  n^{-1} G\trans G & \Pconv H \\
  \that & \Pconv m
\end{align}
\end{subequations}
where $H$ and $m$ are both non-random limits.  Substituting into \eqref{eq:exact} and rearranging,
\begin{equation}\label{eq:flexlim}
  \flex - \frac{ d \log n }{ 2 } \Pconv 
\frac{1}{2} \Big\{ {-d} (\log \sigma^2 + \log \lambda^2) + \log \det H + \lambda^2 \, \Vert m \Vert^2 \Big\} .
\end{equation}
The second term on the lefthand side is the Bayesian Information Criterion
(BIC) penalty.  Thus, for the Gaussian linear model and IID sampling,
\begin{equation}\label{eq:BIC}
  \flex = \text{BIC penalty} + \bigO_p(1) , \qquad \text{$d$ fixed, $n \to \infty$.}
\end{equation}
This large-$n$ calculation
is also applicable for non-Gaussian models which satisfy the conditions 
for asymptotic posterior Normality \citep[see, e.g.,][sec.~7.4]{schervish95}.

However, we advise caution.  First, as already noted,  we cannot
presume an approximately Gaussian posterior distribution in modern
practice, and therefore
\[
\mbox{\flex $\approx$ BIC penalty}
\]
is a poor generic approximation.  Second, it is not safe to drop
$\bigO_p(1)$ terms in model comparison \citep{gelfand94}.  In the
Gaussian linear model, \eqref{eq:flexlim} shows that \flex is
approximately the sum of a term in $\log n$ and a term tending to
a constant.  But $\log n$ grows very slowly: as memorably expressed
by \citet[p.~117]{gowers02}, $\log n$ is about $2.3$~times the number
of digits in the decimal expression of $n$.  It could take a truely
massive $n$ for the $\log n$ term to dominate the \flex in \emph{every
one} of the models in $\curlyM$, as required to be reasonably sure that
\flex and the BIC penalty would give the same outcome.  In practice, choosing to
approximate \flex with the BIC penalty is effectively a choice to use a
different penalty.  In the terms of section~\ref{sec:candidate}, this
is equivalent to penalizing the log-evidence.  Unless the application
provides a specific justification for penalizing the log-evidence,
it is more defensible to estimate the evidence directly.

\section{Summary}
\label{sec:summary}

There are many issues to consider in model selection, starting
with whether it is even sensible to select a single model from a
set of candidates.  Pragmatically, though, operating with a single
data-selected model is a powerful simplification in both inference
and prediction.  Accepting this, the next issue is what loss
function is appropriate for selecting a single model.  While it is
not obligatory to formulate model selection as a decision problem, it is
hardly defensible to make a choice without contemplating the consequences
of an error.  The zero-one loss function used in this paper is a
reasonable default choice, but if the candidate models live in a
metric space there would be alternatives which could use the distances
between models more effectively.

After these two steps, selecting the single model which maximizes the
evidence is `hopefully' the Bayes Rule and definitely admissible,
under zero-one loss.  This is the basis of the Bayesian argument for
maximizing the evidence, and a reassuring feature for all statisticians.

But we have also described maximizing the evidence from a Frequentist
perspective, that of penalizing fit using a measure of model
complexity.  The general justification of this penalization approach
is heuristic---avoiding the dangers of overfitting---and therefore
there cannot be a `right' complexity penalty.
Our contribution in this paper is to clarify that there is a unique
measure of complexity, which we term `flexibility', under which the
log-evidence decomposes exactly into `fit minus flexibility'.
Thus there is a second argument for maximizing the evidence, if
`flexibility' is accepted as a reasonable way to quantify and
penalize model complexity.

These are two arguments to the conclusion that maximizing
the evidence is a sensible way to allow the data to select
a single model, one Bayesian and one Frequentist.  Neither argument
on its own is compelling, but there must be some appeal in a criterion
which has
`cross-party' support.  As we describe in Section~\ref{sec:candidate},
the decision \emph{not} to use flexibility as the complexity penalty
should raise questions in the mind of an auditor which might be better
unasked---unless the statistician has a good answer.  This is a weak
justification, to be sure, but model selection and model complexity are
subtle and ambiguous topics, and in many applications we might welcome
even a weak justification.

Finally, if flexibility is accepted as a reasonable way to quantify
and penalize model complexity, then we strongly recommend
estimating the evidence directly, rather than using \mbox{`evidence
equals fit minus flexibility'} and then replacing flexibility with a
simpler term such as the BIC penalty.

\section*{Acknowledgments}

We began this paper while CEP was Heilbronn Distinguished Visitor in Data
Science at the University of Bristol, UK, in Spring 2019.  We would
like to thank the two reviewers for their detailed comments on previous
versions of this paper, which improved it considerably.

\singlespacing 


\end{document}